\magnification=1200
\def\C{{\bf C}}
\def\N{{\bf N}}
\def\R{{\bf R}}
\def\hal{{\vrule height 10pt width 4pt depth 0pt}}

\centerline{The Kadison-Singer problem}

\centerline{in discrepancy theory}
\medskip

\centerline{Nik Weaver\footnote{*}{Partially supported by NSF grant DMS-0070634
\hfill\break
Math Subject Classification numbers: Primary 05A99, 11K38, 46L05}}
\bigskip
\bigskip

{\narrower{
\noindent \it We give a combinatorial form of
the Kadison-Singer problem, a famous problem in C*-algebra. This
combinatorial problem, which has several minor variations, is a
discrepancy question about vectors in $\C^n$. Some partial
results can be easily deduced from known facts in discrepancy
theory.
\bigskip}}
\bigskip

In its original form, the so-called Kadison-Singer problem [8] asks whether
every pure state on an atomic maximal abelian self-adjoint subalgebra
of $B(H)$ extends uniquely to a pure state on $B(H)$, where $H$ is
a separable Hilbert space and $B(H)$ is the C*-algebra of bounded
linear operators on $H$. It is considered a basic question about
the most fundamental nonabelian C*-algebra, and has generated
a fairly substantial research literature. (See [1], [5], and [11]
for references.)
\medskip

It has been known since [8] that the problem can be reformulated as
a question about finite complex matrices. Let a {\it diagonal projection}
be a matrix whose off-diagonal entries are zero and whose diagonal
entries are each either zero or one. The ``paving problem'' form of
the Kadison-Singer problem asks, for each $\epsilon > 0$, to find a
natural number $r$ such that the following holds: for any complex
$n\times n$ matrix $A$ whose diagonal is zero, there exist $n \times n$
diagonal projections $Q_1, \ldots, Q_r$ which sum to the identity matrix,
such that $\|Q_jAQ_j\| \leq \epsilon \|A\|$ for all $j$. Note that
the number $r$ of projections must be independent of $n$. The norm
used here is the operator norm for matrices acting on Euclidean
space, i.e., $\|A\| = \sup\{\|Av\|_2: \|v\|_2 = 1\}$. A short
proof of the equivalence of this question with the original problem
is given in [10].
\medskip

This version of the problem already has a discrepancy-theoretic
flavor, and partial results on it were obtained by
Bourgain and Tzafriri using probabilistic methods [5]. However,
we base our approach on a different reduction due to Akemann and
Anderson, which resembles the above but in which $A$ becomes
an orthogonal projection with near-zero diagonal.
\medskip

A complex $n \times n$ matrix $P$ is an {\it orthogonal projection}
if $P^2 = P^* = I_n$, where ${}^*$ denotes Hermitian adjoint and
$I_n$ is the $n\times n$ identity matrix. For such a matrix
$P = [p_{ij}]$, let $\delta(P) = \max_i p_{ii}$. Akemann and Anderson
[1] considered the conjecture that there exist $\epsilon,
\delta > 0$ with the following property: for any complex $n \times n$
orthogonal projection $P$ with $\delta(P) \leq \delta$, there is a diagonal
projection $Q$ such that $\|QPQ\| \leq 1 - \epsilon$ and
$\|(I_n - Q)P(I_n - Q)\| \leq 1 - \epsilon$.
(As above, $\|\cdot\|$ denotes operator norm.) They showed that this
conjecture, if true, would imply a positive solution to the Kadison-Singer
problem. Note that $\epsilon$ and $\delta$ must be independent of $n$.
\bigskip
\bigskip

\noindent {\bf 1. Combinatorial versions of the Kadison-Singer problem}
\bigskip

We now state a more directly combinatorial version of the Kadison-Singer
problem.
\bigskip

\noindent {\bf Notation.}  Let $e_1, \ldots, e_k$ be the canonical
orthonormal basis of $\C^k$, let $\|v\|_2$ denote the Euclidean norm
of $v \in \C^k$, let $I_k$ be the $k \times k$ identity matrix, and
for $v \in \C^k$ let $A_v: \C^k \to \C^k$ be the rank one operator
$A_v: u \mapsto \langle u,v\rangle v$. (So $|\langle u,v\rangle|^2
= \langle A_v u,u\rangle$.)
\bigskip

For any natural number $r \geq 2$ we have the following conjecture.
\bigskip

\noindent {\bf Conjecture KS${}_r$.} There exist universal constants $N \geq 2$
and $\epsilon > 0$ such that the following holds. Let $v_1, \ldots, v_n
\in \C^k$ satisfy $\|v_i\|_2 \leq 1$ for all $i$ and suppose
$$\sum_i |\langle u,v_i\rangle|^2 \leq N$$
for every unit vector $u \in \C^k$. Then there exists a partition
$X_1, \ldots, X_r$ of $\{1, \ldots, n\}$ such that
$$\sum_{i \in X_j} |\langle u,v_i\rangle|^2 \leq N - \epsilon$$
for every unit vector $u \in \C^k$ and all $j$.
\bigskip

Note that $N$ and $\epsilon$ must be independent of $n$ and $k$.
\medskip

Also, observe that for any $X \subset \{1, \ldots, n\}$ we have
$0 \leq \sum_X A_{v_i} \leq N\cdot I_k$. Letting $A$ be this sum,
we have $\|A\| \leq N$ and so
$$\left|\sum_X |\langle u, v_i\rangle|^2
- \sum_X |\langle u,' v_i\rangle|^2\right|
= \big|\langle Au, u\rangle - \langle Au', u'\rangle\big|
\leq 2N\|u - u'\|_2$$
for any unit vectors $u$ and $u'$. It follows that the conclusion of
KS${}_r$ really only needs to be verified on an $\epsilon/4N$-net in
the unit sphere of $\C^k$. This remark is due to Nets Katz.
\bigskip

\noindent {\bf Theorem 1.} {\it The Kadison-Singer problem has a
positive solution if and only if Conjecture KS${}_r$ is true for some
$r \geq 2$.}
\medskip

\noindent {\it Proof.} Suppose Conjecture KS${}_r$ holds
for some $r$, $N$, and $\epsilon$. We claim that for any
complex $n \times n$ orthogonal projection $P$ with $\delta(P) \leq
1/N$ (this notation was defined in the introduction)
there exist $n \times n$ diagonal projections $Q_1, \ldots, Q_r$
which sum to the identity and satisfy $\|Q_jPQ_j\| \leq 1 - \epsilon/N$
for all $j$. To see this, let $P$ be a complex $n\times n$ orthogonal
projection with $\delta(P) \leq 1/N$. If $P$ has rank $k$ then its
range is a $k$-dimensional subspace $V \subset \C^n$. Define
$v_i = \sqrt{N}\cdot Pe_i \in V$ for $1 \leq i \leq n$. Observe that
$$\|v_i\|_2^2 = N\cdot\|Pe_i\|_2^2 = N\langle Pe_i, e_i\rangle
\leq N\delta(P) \leq 1$$
for all $i$. Also, for any unit vector $u \in V$ we have
$$\sum_i |\langle u,v_i\rangle|^2
= \sum_i |\langle u,\sqrt{N}Pe_i\rangle|^2
= N\cdot \sum_i |\langle u,e_i\rangle|^2 = N.$$
Thus, Conjecture KS${}_r$ asserts that there exists a partition
$X_1, \ldots, X_r$ of $\{1, \ldots, n\}$ such that
$$\sum_{i \in X_j} |\langle u,v_i\rangle|^2 \leq N - \epsilon$$
for every unit vector $u \in V$ and all $j$. Let $Q_j$ be the
$n\times n$ diagonal projection defined by
$$Q_j e_i = \cases{e_i& if $i \in X_j$\cr 0& if $i \not\in X_j$\cr}$$
($1 \leq j \leq r$, $1 \leq i \leq n$). Then $Q_1 + \cdots + Q_r = I_n$,
and for any unit vector $u \in V$ we have
$$\eqalign{\|Q_jP u\|_2^2 &= \sum_i |\langle Q_jP u, e_i\rangle|^2
= \sum_i |\langle u, PQ_je_i\rangle|^2\cr
&= N^{-1} \sum_{i \in X_j} |\langle u,v_i\rangle|^2
\leq 1 - \epsilon/N.\cr}$$
This shows that $\|Q_jPQ_j\| = \|Q_jP\|^2 \leq 1 - \epsilon/N$ for
all $j$, as claimed.
\medskip

The claim implies a positive solution to the Kadison-Singer problem
by a minor modification of Propositions 7.6 and 7.7 of [1]. (Those
results are stated for the case $r = 2$, but generalize to arbitrary
$r$ with the trivial observation that for any ultrafilter $\cal U$
on $\N$ and any partition $Y_1, \ldots, Y_r$ of $\N$, we have
$Y_j \in \cal U$ for some $j$.) Thus, the reverse implication holds.
\medskip

Conversely, suppose Conjecture KS${}_r$ fails for all $r$. Fix
$N = r \geq 2$ and let $v_1, \ldots, v_n \in \C^k$ be a counterexample
with $\epsilon = 1$. Let $w_i = v_i/\sqrt{N}$ and observe that
$\|A_{w_i}\| = \|w_i\|_2^2 \leq 1/N$ for all $i$ and $\sum_1^n A_{w_i}
\leq I_k$. Then $I_k - \sum A_{w_i}$ is a positive finite rank operator,
so we can find positive rank one operators $A_{w_i}$ ($n + 1 \leq i \leq m$)
such that $\|A_{w_i}\| \leq 1/N$ for all $i$ and $\sum_1^m A_{w_i} = I_k$.
\medskip

Define an embedding $\Phi: \C^k \to \C^m$ by $\langle \Phi u, e_i\rangle
= \langle u, w_i\rangle$ for $1 \leq i \leq m$. For any $u \in \C^k$ we have
$$\|\Phi u\|_2^2 = \sum_1^m |\langle \Phi u, e_i\rangle|^2
= \sum_1^m |\langle u,w_i\rangle|^2
= \sum_1^m \langle A_{w_i} u, u\rangle = \|u\|_2^2,$$
so $\Phi$ is isometric. Let $P$ be the orthogonal projection in
$\C^{n+m}$ with range $\Phi(\C^k)$; then
$$\langle Pe_i, \Phi w_j\rangle = \langle e_i, \Phi w_j\rangle
= \langle w_i, w_j\rangle = \langle \Phi w_i, \Phi w_j\rangle$$
for all $i$ and $j$, which shows that $Pe_i = \Phi w_i$ since the
$w_i$ clearly span $\C^k$. Let $D$ be
the diagonal matrix with the same diagonal as $P$, i.e., $d_{ii} = p_{ii}$
($1 \leq i \leq m$). Then $\|D\| = \max_i \|w_i\|_2^2 \leq 1/N$.
\medskip

Let $Q_1, \ldots, Q_r$ be any $m \times m$ diagonal projections
which sum to the identity matrix. Define a partition $X_1, \ldots, X_r$
of $\{1, \ldots, m\}$ by letting $X_j$ be the diagonal of $Q_j$.
According to our choice of
$v_1, \ldots, v_n$, we infer that there exists $1 \leq j \leq r$ and
$u \in \C^k$, $\|u\|_2 = 1$, such that $\sum_{i \in X_j \cap \{1, \ldots, n\}}
|\langle u, v_i\rangle|^2 > N - 1$, so that
$\sum_{i \in X_j} |\langle u, w_i\rangle|^2 > 1 - 1/N$. It follows that
$$\eqalign{\|Q_j P(\Phi u)\|_2^2
&\geq \sum_{i = 1}^m |\langle Q_j P(\Phi u), e_i\rangle|^2
= \sum_{i \in X_j} |\langle \Phi u, e_i\rangle|^2\cr
&= \sum_{i \in X_j} |\langle u, w_i\rangle|^2
> 1 - 1/N\cr}$$
for this value of $j$. Thus $\|Q_jPQ_j\| = \|Q_jP\|^2 > 1 - 1/N$. Finally,
the matrix $A = P - D$ has zero diagonal and satisfies $\|A\| \leq 1 + 1/N$,
and the preceding shows that for any $m \times m$ diagonal
projections $Q_1, \ldots, Q_r$ which sum to the identity we have
$$\|Q_jAQ_j\| \geq \|Q_jPQ_j\| - \|Q_jDQ_j\|
\geq 1 - 2/N$$
for some $j$.
Thus, as $N  = r \to \infty$, we obtain a sequence of examples which
falsify the paving problem version of the Kadison-Singer problem given
in the introduction.\hfill\hal
\bigskip

Now we indicate possible modifications in Conjecture KS${}_r$ which do
not alter its truth-value.
\bigskip

\noindent {\bf Theorem 2.} {\it If either or both of the following
modifications is made to Conjecture KS${}_r$, the resulting
conjecture is equivalent to Conjecture KS${}_r$:
\medskip

\narrower{
\noindent (a) require $\epsilon = 1$;
\medskip

\noindent (b) assume $\sum_i |\langle u,v_i\rangle|^2 = N$ for every
unit vector $u$ instead of $\sum_i |\langle u,v_i\rangle|^2 \leq N$.
\medskip}}

\noindent {\it Proof.} If Conjecture KS${}_r$ holds for $\epsilon = 1$
then it obviously holds for some $\epsilon > 0$. Conversely, suppose
it holds for some $N$ and $\epsilon$. Since it remains true for all
smaller $\epsilon$, we may assume $1/\epsilon \geq 1$ is an integer.
Then scaling the vectors $v_i$ by $\sqrt{\epsilon}$ shows that it
remains true with $N/\epsilon$ in place of $N$ and $1$ in place of
$\epsilon$ (and $\|v_i\|_2 \leq 1/\sqrt{\epsilon}$ for all $i$, so
certainly for $\|v_i\|_2 \leq 1$). Thus, mandating $\epsilon = 1$
does not change the truth of the conjecture.
\medskip

Now we consider modification (b). We will show that the truth of
Conjecture KS${}_r$ for some $N$ and $\epsilon$ is equivalent to its
truth with this modification, for the same values of $N$ and $\epsilon$.
By the last paragraph, it follows that also including modification (a)
has no effect.
\medskip

Conjecture KS${}_r$ trivially implies the conjecture with modification
(b) for the same values of $N$ and $\epsilon$. Conversely, suppose the
conjecture holds with modification (b), for some given values of $N$ and
$\epsilon$. Let $v_1, \ldots, v_n \in \C^k$ satisfy $\|v_i\|_2 \leq 1$ for
all $i$ and suppose $\sum_i |\langle u,v_i\rangle|^2 \leq N$ for every unit
vector $u \in \C^k$. Then $\sum_i A_{v_i} \leq N\cdot I_k$, so the operator
$N\cdot I_k - \sum_i A_{v_i}$ is a positive finite rank operator, and just
as in the proof of Theorem 1 we can find positive rank one operators
$A_{v_i}$ ($n + 1 \leq i \leq m$) such that $\|A_{v_i}\| \leq 1$ for all
$i$ and $\sum_1^m A_{v_i} = N\cdot I_k$. Then the vectors $v_1, \ldots, v_m$
satisfy the modified hypotheses of the conjecture, so we infer the existence
of a partition $X_1, \ldots, X_r$ of $\{1, \ldots, m\}$ such that
$\sum_{X_j} |\langle u, v_i\rangle|^2 \leq N - \epsilon$ for every unit
vector $u \in \C^k$ and all $j$. Letting $Y_j = X_j \cap \{1, \ldots, n\}$,
we obtain $\sum_{Y_j} |\langle u,v_i\rangle|^2 \leq N - \epsilon$ for every
unit vector $u \in \C^k$ and all $j$. We conclude that Conjecture KS${}_r$
holds for the original vectors $v_1, \ldots, v_n$.\hfill\hal
\bigskip

Conjecture KS${}_r$ can also be modified so that the vectors $v_i$ must
have unit length, though at a significant cost to $\epsilon$.
\bigskip

\noindent {\bf Conjecture KS${}_r{}'$.} There exist universal constants
$N \geq 4$ and $\epsilon > \sqrt{N}$ such that the following holds. Let
$v_1, \ldots, v_n \in \C^k$ satisfy $\|v_i\|_2 = 1$ for all $i$ and suppose
$$\sum_i |\langle u,v_i\rangle|^2 \leq N$$
for every unit vector $u \in \C^k$. Then there exists a partition
$X_1, \ldots, X_r$ of $\{1, \ldots, n\}$ such that
$$\sum_{i \in X_j} |\langle u,v_i\rangle|^2 \leq N - \epsilon$$
for every unit vector $u \in \C^k$ and all $j$.
\bigskip      

\noindent {\bf Theorem 3.} {\it The following are equivalent:
a positive solution of the Kadison-Singer problem; the truth of
Conjecture KS${}_r{}'$ for some $r \geq 2$; and the truth of
Conjecture KS${}_r{}'$ with modification (b) of Theorem 2, for
some $r \geq 2$.}
\medskip

\noindent {\it Proof.} Suppose KS${}_r{}'$ fails for all $r$. Then
we can construct a sequence of matrices $A$ which cannot be paved exactly
as in the proof of Theorem 1. In this case we obtain that $\|A\| \leq 1 + 1/N$
and for any diagonal projections $Q_1, \ldots, Q_r$ which sum to the identity
we have $\|Q_jAQ_j\| > 1 - 1/\sqrt{N} - 2/N$ for some $j$. This is not as
sharp as the conclusion in Theorem 1, but it is sufficient to establish a
negative solution to the paving problem form of the Kadison-Singer problem.
Thus a positive solution of the Kadison-Singer problem implies the truth
of Conjecture KS${}_r{}'$ for some $r \geq 2$.
\medskip

Next, the truth of Conjecture KS${}_r{}'$ without modification (b) clearly
implies its truth with modification (b).
\medskip

Finally, assume that Conjecture KS${}_r{}'$ holds with modification (b), for
some $r \geq 2$, $N \geq 4$, and $\epsilon > \sqrt{N}$. We will verify that
Conjecture KS${}_r$ holds for $r$, $N - \sqrt{N}$, and $\epsilon - \sqrt{N}$;
this is sufficient. (Although $N - \sqrt{N}$ need not be an integer, this is
not a problem; a scaling argument as in the first part of the proof of
Theorem 2 can then be used to establish the truth of KS${}_r$ for all
integers larger than $N - \sqrt{N}$.)
\medskip

Let $v_1, \ldots, v_n \in \C^k$ satisfy $\|v_i\|_2 \leq 1$ for all $i$
and suppose $\sum_i |\langle u,v_i\rangle|^2 \leq N - \sqrt{N}$ for every
unit vector $u \in \C^k$. Replacing $\C^k$ with ${\rm span}\{v_i\}$ if
necessary, we may assume $n \geq k$. Now let $m = k + n$ and define unit
vectors $w_i$ ($1 \leq i \leq n$) in $\C^m \cong \C^k \oplus \C^n$ by
$w_i = v_i + \sqrt{1 - \|v_i\|_2^2}\, e_{k + i}$.
It is clear that $\|w_i\|_2 = 1$ for all $i$.
Also define $w_{n+i}$ for $1 \leq i \leq k$ by $w_{n+i} = e_{k+i}$.
\medskip

Let $P_1$ and $P_2$ be the orthogonal projections of $\C^m$ onto
$\C^k$ and $\C^n$ in the decomposition $\C^m \cong \C^k \oplus \C^n$.
For any unit vector $u \in \C^m$ we have
$$\eqalign{\sum_1^m |\langle u, w_i\rangle|^2
&= \sum_1^n |\langle u, v_i\rangle|^2
+ \sum_1^n (1 - \|v_i\|_2^2)|\langle u, e_{k+i}\rangle|^2\cr
&{}\qquad\qquad + 2{\rm Re}\sum_1^n \sqrt{1 - \|v_i\|_2^2}
\langle u, v_i\rangle\overline{\langle u, e_{k+i}\rangle}
+ \sum_1^k |\langle u, e_{k+i}\rangle|^2.\cr}$$
Now
$$\eqalign{\sum_1^n |\langle u, v_i\rangle|^2
&+ \sum_1^n (1 - \|v_i\|_2^2)|\langle u, e_{k+i}\rangle|^2
+ \sum_1^k |\langle u, e_{k+i}\rangle|^2\cr
& \leq \sum_1^n |\langle P_1u, v_i\rangle|^2
+ 2\sum_1^n |\langle P_2u, e_{k+i}\rangle|^2\cr
& \leq (N - \sqrt{N})\|P_1u\|_2^2 + 2\|P_2u\|_2^2 \leq N - \sqrt{N}\cr}$$
since $N \geq 4$, and by the Cauchy-Schwarz inequality
$$\eqalign{&2{\rm Re}\sum_1^n \sqrt{1 - \|v_i\|_2^2}
\langle u, v_i\rangle\overline{\langle u, e_{k+i}\rangle}\cr
&\qquad\qquad\leq 2\left(\sum_1^n |\langle u, v_i\rangle|^2\right)^{1/2}
\left(\sum_1^n (1 - \|v_i\|_2^2)|\langle u, e_{k+i}\rangle|^2\right)^{1/2}\cr
&\qquad\qquad\leq 2\sqrt{N - \sqrt{N}}\|P_1u\|_2\|P_2u\|_2
\leq \sqrt{N - \sqrt{N}}.\cr}$$
Thus, we conclude that $\sum_1^m |\langle u, w_i\rangle|^2 \leq
N - \sqrt{N} + \sqrt{N - \sqrt{N}} \leq N$.
\medskip

Let $B = N \cdot I_m - \sum_1^m A_{w_i}$. The preceding shows that
$B$ is a positive operator,
and ${\rm tr}(B) = (N-1)m$ since ${\rm tr}(A_{w_i}) = \|w_i\|_2^2 = 1$ for
all $i$. Let $\{f_t\}$ be an orthonormal basis of $\C^m$ which diagonalizes
$B$ and say $Bf_t = b_tf_t$ ($1 \leq t \leq m$). Notice that
$\sum_1^m b_t = {\rm tr}(B) = (N-1)m$. Then define unit vectors
$u_1, \ldots, u_m \in \C^m$ by $\langle u_s, f_t\rangle =
\sqrt{b_t/(N-1)m}\, e^{2\pi ist}$. We have
$$\sum_s \langle f_t, u_s\rangle u_s = {{b_t}\over{N-1}}f_t$$
for all $1 \leq t \leq m$, that is, $\sum_s A_{u_s} = B/(N-1)$.
Thus, letting $w_{m+1}, \ldots, w_{mN}$ consist of $N-1$
copies of each of the vectors $u_1, \ldots, u_m$, we have
$\|w_i\|_2 = 1$ for $1 \leq i \leq mN$ and
$$\sum_1^{mN} A_{w_i} = N\cdot I_m.$$
We are now in a position to apply Conjecture KS${}_r{}'$ with modification
(b). We infer that there exists a partition $X_1, \ldots, X_r$ of
$\{1, \ldots, mN\}$ such that
$$\sum_{i \in X_j} A_{w_i} \leq (N - \epsilon)I_m$$
for all $j$. Letting $Y_j = X_j \cap \{1, \ldots, n\}$, we obtain
$$\sum_{i \in Y_j} A_{v_i} \leq (N - \epsilon)I_k
= ((N - \sqrt{N}) - (\epsilon - \sqrt{N}))I_k$$
for all $j$. We conclude that Conjecture KS${}_r$ holds for $r$,
$N - \sqrt{N}$, and $\epsilon - \sqrt{N}$, as desired.\hfill\hal
\bigskip

It is unclear whether the real version of Conjecture KS${}_r$ or any of its
variants is substantially different from the complex version. These variants
remain equivalent by the same proofs in the real case, with the one possible
exception of Conjecture KS${}_r{}'$ with modification (b), whose preceding
equivalence proof does use complex scalars.
\bigskip
\bigskip

\noindent {\bf 2. Positive results}
\bigskip

We single out the case $r = 2$ for special attention:
\bigskip

\noindent {\bf Conjecture KS${}_2$.} There exist universal constants $N \geq 2$
and $\epsilon > 0$ such that the following holds. Let $v_1, \ldots, v_n
\in \C^k$ satisfy $\|v_i\|_2 \leq 1$ for all $i$ and suppose
$$\sum_i |\langle u,v_i\rangle|^2 \leq N$$
for every unit vector $u \in \C^k$.
Then for some choice of signs we have
$$\left|\sum_i \pm |\langle u,v_i\rangle|^2\right|
\leq N - \epsilon$$
for every unit vector $u \in \C^k$.
\bigskip

Again, $N$ and $\epsilon$ must be independent of $n$ and $k$.
\medskip

We have altered the statement slightly to conform more closely in style
to traditional discrepancy statements. Equivalence to the $r = 2$ case of
Conjecture KS${}_r$ as stated above is an easy exercise. (The
value of $\epsilon$ changes by a factor of 2.) Conjecture KS${}_2$ is
equivalent to Conjecture 7.1.3 of [1].
\medskip

In this section we present three positive partial results on
Conjecture KS${}_2$ which follow from known general results. For
background on discrepancy theory, see [6] or [9].
\medskip

First we observe that a strong form of the conclusion of Conjecture
KS${}_2$ always holds on an orthonormal basis, if not for all unit
vectors in $\C^k$. Note that the hypothesis $\sum |\langle u, v_i\rangle|^2
\leq N$ is not needed for this.
\bigskip

\noindent {\bf Proposition 4.} {\it Let $v_1, \ldots, v_n \in \C^k$ satisfy
$\|v_i\|_2 \leq 1$ for all $i$. Then there is a choice of signs such that
$$\left| \sum_i \pm |\langle u,v_i\rangle|^2\right| \leq 2$$
for all $u \in \{e_1, \ldots, e_k\}$.}
\medskip

\noindent {\it Proof.} We have $\sum_j |\langle e_j,v_i\rangle|^2
= \|v_i\|_2^2 \leq 1$
for all $i$. Thus the vectors $a_1, \ldots, a_n \in \R^k$ defined by
$\langle a_i, e_j\rangle = |\langle e_j,v_i\rangle|^2$ satisfy
$\|a_i\|_1 \leq 1$
for all $i$. It follows from the Beck-Fiala theorem [4]
that there is a choice of signs such that
$$\left\| \sum_i \pm a_i\right\|_\infty \leq 2,$$
i.e., $|\sum \pm \langle a_i, e_j\rangle| \leq 2$ for all $j$.
Since $\sum_i \pm |\langle e_j, v_i\rangle|^2 =
\sum_i \pm \langle a_i, e_j\rangle$, we are done.\hfill\hal
\bigskip

Next, we show that the conclusion of Conjecture KS${}_2$ can be achieved
for arbitrary $N$ if $\epsilon$ is allowed to depend on $n$. The result
is surprisingly difficult; we prove it using a clever theorem on matroid
partitions. The same theorem was used to a similar purpose in [1].
\bigskip

\noindent {\bf Proposition 5.} {\it Let $N \geq 2$, let
$v_1, \ldots, v_n \in \C^k$ satisfy $\|v_i\|_2 \leq 1$ for all $i$,
and suppose
$$\sum_i |\langle u,v_i\rangle|^2 \leq N$$
for every unit vector $u \in \C^k$.
Then there is a choice of signs such that
$$\left| \sum_i \pm |\langle u,v_i\rangle|^2\right| < N$$
for every unit vector $u \in \C^k$.}
\medskip

\noindent {\it Proof.} The proof requires the stronger hypothesis
$\sum_i |\langle u,v_i\rangle|^2 = N$ for every unit vector $u$.
As in the proofs of Theorems
1 and 2, we can achieve this hypothesis by enlarging the set of vectors.
Let $v_1, \ldots, v_m$ ($m \geq n$) be an expanded list which satisfies
$\|v_i\|_2 \leq 1$ for all $i$ and $\sum_1^m |\langle u,v_i\rangle|^2 = N$
for every unit vector $u \in \C^k$.
\medskip

The collection of subsets of $\{v_1, \ldots, v_m\}$ which are linearly
independent in $\C^k$ constitutes a matroid. (See [7] for definitions.)
Now for any subset
$X \subset \{1, \ldots, m\}$ with cardinality $|X|$, let $V =
{\rm span}\{v_i: i \not\in X\}$ and let $d$ be the dimension of $V$. Then
$$\sum_{i \not\in X} A_{v_i} \leq N\cdot P_V$$
where $P_V$ is the orthogonal projection onto $V$, so
$$\sum_{i \not\in X} {\rm tr}(A_{v_i}) \leq Nd.$$
Since $\sum_1^m {\rm tr}(A_{v_i}) = {\rm tr}(N\cdot I_k) = Nk$, combining
the preceding with
$$\sum_{i \in X} {\rm tr}(A_{v_i}) \leq |X|$$
(since ${\rm tr}(A_{v_i}) = \|v_i\|_2^2 \leq 1$ for all $i$) yields
$N(k - d) \leq |X|$. This verifies the hypothesis of the Edmonds-Fulkerson
theorem ([7], Theorem 2c), and we deduce that $\{v_1, \ldots, v_m\}$ can
be partitioned into two sets (in fact $N$ sets, but this does not seem to
help matters any) $X_1$
and $X_2$, each of which spans $\C^k$. It follows that the quantity
$\sum_{i \in X_j} |\langle u,v_i\rangle|^2$ is never zero as $u$ ranges
over all unit vectors in $\C^k$ ($j = 1,2$), and therefore
$\sum_{i \in X_j} |\langle u,v_i\rangle|^2 < N$ for every unit vector
$u \in \C^k$ ($j = 1,2$). The same final conclusion obviously holds for
sums over $Y_j = X_j \cap \{1, \ldots, n\}$ ($j = 1,2$), which is
enough.\hfill\hal
\bigskip

By compactness, the conclusion of Proposition 5 can be strengthened to
say that $|\sum_i \pm |\langle u, v_i\rangle|^2| \leq N - \epsilon$
for some $\epsilon > 0$ and every unit vector $u$. The existence of a
universal $\epsilon$ for any fixed value of $n$ then follows by another
easy compactness argument. (Recall that we can assume $k \leq n$.)
\medskip

Lastly, we observe that for fixed $k$ the conclusion of conjecture KS${}_2$
can be achieved with no assumption on $\sum |\langle u, v_i\rangle|^2$.
\bigskip

\noindent {\bf Proposition 6.} {\it For fixed $k$, there exists $M = M(k)$
such that the following holds.  Let $v_1, \ldots, v_n \in \C^k$
satisfy $\|v_i\|_2 \leq 1$ for all $i$. Then there is a choice of signs
such that
$$\left| \sum_i \pm |\langle u,v_i\rangle|^2\right| \leq M$$
for every unit vector $u \in \C^k$.}
\medskip

\noindent {\it Proof.} We use a vector balancing theorem due to
Banaszczyk [3]. Let $M_k^{sa}(\C)$ be the real vector space of
self-adjoint $k\times k$ matrics, with Euclidean (Hilbert-Schmidt)
norm given by $\|A\|_2 = ({\rm tr}(A^*A))^{1/2}$. Define $R > 0$ by
the condition that the set of matrices in $M_k^{sa}(\C)$ with operator
norm at most $R$ has Gaussian measure $1/2$. Let $K$ be this set and let
$M = 5R$. Also let $B_i = {1\over 5}A_{v_i}$ for all $i$; then
$\|B_i\|_2 \leq 1/5$ for all $i$ and Banaszczyk's theorem asserts
that there is a choice of signs such that $\sum \pm B_i \in K$.
It follows that
$$\left\|\sum \pm A_{v_i}\right\| \leq M,$$
and hence
$$-M\cdot I_k \leq \sum \pm A_{v_i}\leq M\cdot I_k,$$
so that $|\sum_i \pm |\langle u,v_i\rangle|^2| =
|\sum_i \pm \langle A_{v_i}u,u\rangle| \leq M$ for every unit vector
$u \in \C^k$.\hfill\hal
\bigskip

Oddly, none of the results used in this section rely on probabilistic
methods. It seems likely that such methods could be used to obtain
further results; see [2, Chapter 12]. However, we have been unable
to do this.
\bigskip
\bigskip

\noindent {\bf 3. A counterexample}
\bigskip

It is interesting to note that neither Proposition 4 nor Proposition 6
requires the assumption $\sum |\langle u,v_i\rangle|^2 \leq N$. This
suggests replacing KS${}_2$ with an even stronger conjecture which
dispenses with this assumption, and essentially
this was done in [1, Conjecture 7.1]. However,
that version of the conjecture was falsified in [11]. Here we present a
sharper version of the counterexample which provides better asymptotics.
\bigskip

\noindent {\bf Example 7.} Let $k \geq 5$ be an integer. Define vectors
$v_1', \ldots, v_{k-1}' \in \C^k$ by
$$v_i' = (k-2)\alpha e_i -
\left(\sum_{{1 \leq j \leq k-1}\atop{j\neq i}} \alpha e_j\right)
+ \beta e_k$$
where $\alpha = (k-1)^{-3/2}$ and $\beta = (k-1)^{-1/2}$. Also define
$v_i = v_i'/\sqrt{\delta}$ ($1 \leq i \leq k - 1$) where
$\delta = \|v_i'\|_2^2 = (2k-3)/(k-1)^2 \sim 2/k$.
\medskip

We claim that $\sum |\langle u,v_i\rangle|^2 \leq N \equiv 1/\delta$
for every unit vector $u \in \C^k$. To see this, let $u \in \C^k$
be a unit vector and write $u = \sum_1^k a_ie_i$. Then
$$\eqalign{\sum_{i=1}^{k-1} |\langle u,v_i'\rangle|^2
&= \sum_i \left|(k-2)\alpha a_i
- \sum_{j\neq i} \alpha a_j + \beta a_k\right|^2\cr
&= \sum_i \left((k-2)^2\alpha^2 + (k-2)\alpha^2\right)|a_i|^2\cr
&{}\qquad\qquad + 2{\rm Re}\sum_{i\neq j}\left(-2(k-2)\alpha^2
+ (k-3)\alpha^2\right)a_i\bar{a}_j\cr
&{}\qquad\qquad + 2{\rm Re}\sum_i\left((k-2)\alpha\beta -
(k-2)\alpha\beta\right)a_i\bar{a}_k
+ (k-1)\beta^2|a_k|^2\cr
&= \sum_i {{k-2}\over{(k-1)^2}}|a_i|^2
- 2{\rm Re}\sum_{i\neq j} {{1}\over{(k-1)^2}}a_i\bar{a}_j + |a_k|^2\cr
&\leq \sum_i {{k-2}\over{(k-1)^2}}|a_i|^2
+ 2\sum_i {{k}\over{(k-1)^2}} |a_i|^2 + |a_k|^2\cr
&\leq 1\cr}$$
(since $k \geq 5$). All sums in this computation have limits 1 and $k-1$.
With $v_i$ in place of $v_i'$ the sum is bounded by $1/\delta$,
so the claim is proven.
\medskip

Now let $X$ be any subset of $\{1, \ldots, k-1\}$, let $X^c$ be its
complement, and let $c = |X|$.
Since $\sum_i A_{v_i'}(e_k) = e_k$, we have
$$\left\| \sum_X A_{v_i'}(e_k) - \sum_{X^c} A_{v_i'}(e_k)\right\|_2
= 2\left\| \sum_X A_{v_i'}(e_k) - {1\over 2}e_k\right\|_2.$$
The $j$th component of $\sum_X A_{v_i'}(e_k)$ is
$\sum_X\langle e_k, v_i'\rangle\langle v_i', e_j\rangle$; we have
$\sum_X |\langle e_k, v_i'\rangle|^2 = c/(k-1)$ and for $1 \leq j \leq k-1$
$$\sum_X\langle e_k, v_i'\rangle\langle v_i', e_j\rangle
= \cases{-c\alpha\beta&if $j \not\in X$\cr
(k-1-c)\alpha\beta&if $j \in X$.\cr}$$
Thus we can estimate
$$\eqalign{\left\| \sum_X A_{v_i'}(e_k) - {1\over 2}e_k\right\|_2^2
& = (k-1-c)c^2\alpha^2\beta^2 + c(k-1-c)^2\alpha^2\beta^2
+ (c/(k-1) - 1/2)^2\cr
&= c(k-1-c)(k-1)\alpha^2\beta^2 + (c/(k-1) - 1/2)^2\cr
&= c(k-1-c)/(k-1)^3 + (c/(k-1) - 1/2)^2.\cr}$$
This value is minimized at $c = (k-1)/2$, which yields
$$\left\| \sum_X A_{v_i'}(e_k) - {1\over 2}e_k\right\|_2
\geq {1\over{2\sqrt{k-1}}}.$$
We conclude that for any choice of signs we have
$$\left\| \sum_i \pm A_{v_i}\right\| \geq
{1\over{\delta\sqrt{k-1}}} \sim \sqrt{k}/2 \sim \sqrt{N/2},$$
which implies the same lower bound on $|\sum \pm |\langle u, v_i\rangle|^2|
= |\sum \pm \langle A_{v_i}u,u\rangle|$
for some unit vector $u \in \C^k$ which depends on the choice of signs.
This completes the example.
\bigskip

Thus, contrary to the special cases in Propositions 4 and 6, in general
if there is no restriction on $\sum |\langle u, v_i\rangle|^2$
then the quantity $|\sum \pm |\langle u, v_i\rangle|^2|$ can be
arbitrarily large for all choices of signs and some $u$.
\medskip

The following definition is standard. For symmetric convex sets $U, V$
in $\R^d$ let $\beta(U,V)$ be the smallest value of $R$ such that for
any $u_1, \ldots, u_n \in U$ ($n$ arbitrary) there is a choice of signs
such that $\sum \pm u_i \in RV$. Various classical results take $U$ and
$V$ to be $l^p$ unit balls of $\R^d$ for various values of $p$. For
instance, the vector version of the Beck-Fiala theorem states that if
$U$ is the $l^1$ unit ball and $V$ is the $l^\infty$ unit ball then
$\beta(U,V) \leq 2$.
\medskip

The questions posed in this paper suggest interest in ``noncommutative
discrepancy'' questions where $\R^d$ is replaced with the real
$k^2$-dimensional vector space $M_k^{sa}(\C)$ of self-adjoint complex
$k\times k$ matrics and the $l^p$ unit ball of $\R^d$ is replaced with
the unit ball of $M_k^{sa}(\C)$ for the Schatten $p$-norm defined by
$$\|A\|_p = ({\rm tr}\, |A|^p)^{1/p}$$
(where $|A| = (A^*A)^{1/2}$) for $1 \leq p < \infty$ and
$\|A\|_\infty = \|A\|$ (operator norm). The operators $A_{v_i}$ in
Example 7 then show that the noncommutative analog of the Beck-Fiala
theorem fails:
\bigskip

\noindent {\bf Theorem 8.} {\it Let
$$U = \{A \in M_k^{sa}(\C): {\rm tr}\, |A| \leq 1\}$$
and
$$V = \{A \in M_k^{sa}(\C): \|A\| \leq 1\}.$$
Then $\beta(U,V)= \Omega(\sqrt{k})$.}
\bigskip

Obtaining information on the value of $\beta(U,V)$ for other norms seems
difficult. The case where $U = V$ is the unit ball for
the Hilbert-Schmidt norm $\|\cdot\|_2$ is classical, because this
is a Euclidean norm; we have $\beta(U,U) = k$. But in all other
cases little is obvious besides the simple observation that restricting
the matrices to be diagonal shows that the value of $\beta$
in such a matrix problem is always bounded below by its value in the
corresponding vector problem.
\bigskip
\bigskip

\noindent {\bf Acknowledgement}
\bigskip

The author wishes to thank Charles Akemann, Nets Katz, Ji\v{r}i Matou\v{s}ek,
and John Shareshian for help and advice.

\bigskip
\bigskip

[1] C.\ A.\ Akemann and J.\ Anderson, Lyapunov theorems
for operator algebras, {\it Mem.\ Amer.\ Math.\ Soc.\ \bf 94} (1991).
\medskip

[2] N.\ Alon and J.\ H.\ Spencer, {\it The Probabilistic Method}
(second edition), Wiley-Interscience, New York, 2000.
\medskip

[3] W.\ Banaszczyk, Balancing vectors and Gaussian measures of
$n$-dimensional convex bodies, {\it Random Structures Algorithms
\bf 12} (1998), 351-360.
\medskip

[4] J.\ Beck and T.\ Fiala, ``Integer-making'' theorems, {\it
Discrete Appl.\ Math.\ \bf 3} (1981), 1-8.
\medskip

[5] J.\ Bourgain and L.\ Tzafriri, On a problem of Kadison and
Singer, {\it J.\ Reine Angew.\ Math.\ \bf 420} (1991), 1-43.
\medskip

[6] B.\ Chazelle, {\it The Discrepancy Method: Randomness and Complexity},
Cambridge University Press, Cambridge, 2000.
\medskip

[7] J.\ Edmonds and D.\ R.\ Fulkerson, Transversals and matroid
partition, {\it J.\ Res.\ Nat.\ Bur.\ Standards Sect.\ B \bf 69B}
(1965), 147-153.
\medskip

[8] R.\ Kadison and I.\ Singer, Extensions of pure states,
{\it Amer.\ J.\ Math.\ \bf 81} (1959), 547-564.
\medskip

[9] J.\ Matou\v{s}ek, {\it Geometric Discrepancy : an Illustrated Guide},
Springer, New York, 1999.
\medskip

[10] B.\ Tanbay, Pure state extensions and compressibility of the
$l_1$-algebra, {\it Proc.\ Amer.\ Math.\ Soc.\ \bf 113}
(1991), 707-713.
\medskip

[11] N.\ Weaver, A counterexample to a conjecture of Akemann and
Anderson, {\it Bull.\ London Math.\ Soc.\ \bf 34} (2002), 1-7.
\bigskip
\bigskip

\noindent Math Dept.

\noindent Washington University

\noindent St.\ Louis, MO 63130 USA

\noindent nweaver@math.wustl.edu

\end